\documentclass[12 pt]{amsart}

\newtheorem{theorem}{Theorem}[section]
\newtheorem{lemma}[theorem]{Lemma}
\newtheorem{proposition}[theorem]{Proposition}
\newtheorem{corollary}[theorem]{Corollary}

\theoremstyle{definition}

\newtheorem{remark}[theorem]{Remark}
\newtheorem{algorithm}[theorem]{Algorithm}

\theoremstyle{remark}

\usepackage{amscd,amssymb}

\begin{document}
\baselineskip 15 pt

\title[Tame Automorphisms Fixing a Variable]
{Tame Automorphisms Fixing a Variable of
Free Associative Algebras of Rank Three}

\author[Vesselin Drensky and Jie-Tai Yu]
{Vesselin Drensky and Jie-Tai Yu}
\address{Institute of Mathematics and Informatics,
Bulgarian Academy of Sciences,
1113 Sofia, Bulgaria}
\email{drensky@math.bas.bg}
\address{Department of Mathematics, the University of Hong Kong,
Hong Kong SAR, China}
\email{yujt@hkucc.hku.hk}

\thanks
{The research of Vesselin Drensky was partially supported by Grant
MI-1503/2005 of the Bulgarian National Science Fund.}

\thanks{The research of Jie-Tai Yu was partially
supported by a Hong Kong RGC-CERG Grant.}

\subjclass[2000]
{Primary 16S10. Secondary 16W20; 16Z05.}
\keywords{Automorphisms of free
algebras, tame automorphisms, tame coordinates,
primitive elements in free algebras}

\maketitle

\centerline{Devoted to the 80th anniversary of Professor Boris Plotkin.}

\begin{abstract} We study automorphisms of
the free associative algebra $K\langle x,y,z\rangle$ over a field $K$
which fix the variable $z$. We describe
the structure of the group of $z$-tame automorphisms and derive algorithms which recognize
$z$-tame automorphisms and $z$-tame coordinates.
\end{abstract}

\section*{Introduction}

Let $K$ be an arbitrary field of any characteristic
and let $K[x_1,\ldots,x_n]$ and $K\langle x_1,\ldots,x_n\rangle$ be,
respectively, the polynomial algebra in $n$ variables
and of the free associative algebra of rank $n$,
freely generated by $x_1,\ldots,x_n$.
We may think of $K\langle x_1,\ldots,x_n\rangle$
as the algebra of polynomials in $n$ noncommuting variables.
The automorphism groups $\text{Aut }K[x_1,\ldots,x_n]$ and $\text{Aut }K\langle x_1,\ldots,x_n\rangle$
are well understood for $n\leq 2$ only. The description is trivial
for $n=1$, when the automorphisms $\varphi$ are defined by
$\varphi(x_1)=\alpha x_1+\beta$, where $\alpha\in K^{\ast}=K\backslash 0$ and
$\beta\in K$. The classical results of Jung--van der Kulk \cite{J, K} for $K[x_1,x_2]$
and of Czerniakiewicz--Makar-Limanov \cite{Cz, ML1, ML2}
give that all automorphisms of $K[x_1,x_2]$ and $K\langle x_1,x_2\rangle$ are tame.
Writing the automorphisms of $K[x_1,\ldots,x_n]$
and $\text{Aut }K\langle x_1,\ldots,x_n\rangle$ as $n$-tuples
of the images of the variables, and using $x,y$ instead of $x_1,x_2$,
this means that $\text{Aut }K[x,y]$ and $\text{Aut }K\langle x,y\rangle$
are generated by the affine automorphisms
\[
\psi=(\alpha_{11}x+\alpha_{21}y+\beta_1,\alpha_{12}x+\alpha_{22}y+\beta_2),
\quad \alpha_{ij},\beta_j\in K,
\]
(and $\psi_1=(\alpha_{11}x+\alpha_{21}y,\alpha_{12}x+\alpha_{22}y)$, the
linear part of $\psi$, is invertible) and the triangular
automorphisms
\[
\rho=(\alpha_1x+p_1(y),\alpha_2y+\beta_2), \quad
\alpha_1,\alpha_2\in K^{\ast},p_1(y)\in K[y],\beta_2\in K.
\]
It turns out that the groups
$\text{\rm Aut }K\langle x,y\rangle$ and $\text{\rm Aut }K[x,y]$ are naturally isomorphic.
As abstract groups they are
described as the free product $A\ast_CB$ of the group $A$
of the affine automorphisms and the group $B$ of triangular automorphisms
amalgamating their intersection $C=A\cap B$.
Every automorphism $\varphi$ of $K[x,y]$ and $K\langle x,y\rangle$
can be presented as a product
\begin{equation}\label{canonical form of automorphisms}
\varphi=\psi_m^{\varepsilon_m}\rho_m\psi_{m-1}\cdots\rho_2\psi_1\rho_1^{\varepsilon_1},
\end{equation}
where $\psi_i\in A$, $\rho_i\in B$ ($\varepsilon_1$ and $\varepsilon_m$ are
equal to 0 or 1), and, if $\varphi$ does not belong to the union of $A$ and $B$,
we may assume that $\psi_i\in A\backslash B$, $\rho_i\in B\backslash A$.
The freedom of the product means that if $\varphi$ has a nontrivial presentation
of this form, then it is different from the identity automorphism.

In the case of arbitrary $n$, the tame automorphisms
are defined similarly, as compositions
of affine and triangular automorphisms. One studies not only the
automorphisms but also the coordinates, i.e.,
the automorphic images of $x_1$.

We shall mention few facts related with the topic of the present paper,
for $z$-automorphisms of $K[x,y,z]$
and $K\langle x,y,z\rangle$, i.e., automorphisms fixing the variable $z$.
For more details we refer to the books by van den Essen \cite{E2},
Mikhalev, Shpilrain, and Yu \cite{MSY}, and our survey article \cite{DY1}.

Nagata \cite{N} constructed the automorphism of $K[x,y,z]$
\[
\nu=(x-2(y^2+xz)y-(y^2+xz)^2z,y+(y^2+xz)z,z)
\]
which fixes $z$. He showed that $\nu$ is nontame, or wild, considered as an automorphism
of $K[z][x,y]$, and conjectured that it is wild
also as an element of $\text{Aut }K[x,y,z]$. This was the beginning of the study of $z$-automorphisms.

It is relatively easy to see (and to decide algorithmically) whether
an endomorphism of $K[z][x,y]$ is an automorphism and whether this automorphism is
$z$-tame, or tame as an automorphism of $K[z][x,y]$.
When $\text{char }K=0$, Drensky and Yu \cite{DY2} presented a simple
algorithm which decides whether a polynomial $f(x,y,z)\in K[x,y,z]$ is a $z$-coordinate
and whether this coordinate is $z$-tame. This provided
many new wild automorphisms and wild coordinates of $K[z][x,y]$.
These results in \cite{DY2} are based on a similar algorithm of Shpilrain and Yu \cite{SY1}
which recognizes the coordinates of $K[x,y]$.
Shestakov and Umirbaev \cite{SU1, SU2, SU3}
established that the Nagata automorphism is wild.
They also showed that every wild automorphism of $K[z][x,y]$ is wild as an
automorphism of $K[x,y,z]$.
Umirbaev and Yu \cite {UY} proved that
the $z$-wild coordinates in $K[z][x,y]$ are
wild also in $K[x,y,z]$.
In this way, all $z$-wild examples in \cite{DY2} give automatically
wild examples in $K[x,y,z]$.

Going to free algebras,
the most popular candidate for a wild automorphism of
$K\langle x,y,z\rangle$ is the example of Anick
$(x+(y(xy-yz),y,z+(zy-yz)y)\in \text{Aut }K\langle x,y,z\rangle$,
see the book by Cohn \cite{C2}, p.~343. It fixes
one variable and its abelianization is a tame automorphism of $K[x,y,z]$.
Exchanging the places of $y$ and $z$, we obtain the automorphism
$(x+z(xz-zy),y+(xz-zy)z,z)$ which fixes $z$ (or a $z$-automorphism),
and refer to it as the Anick automorphism.
It is linear in $x$ and $y$, considering $z$ as a ``noncommutative constant''.
Drensky and Yu \cite{DY3} showed that such $z$-automorphisms are $z$-wild
if and only if a suitable invertible $2\times 2$ matrix with entries from $K[z_1,z_2]$
is not a product of elementary matrices. In particular, this gives that the Anick
automorphism is $z$-wild. When $\text{char }K=0$,
Umirbaev \cite{U2} described
the defining relations of the group of tame automorphisms of $K[x,y,z]$.
He showed that $\varphi=(f,g,h)\in \text{Aut }K\langle x,y,z\rangle$ is wild
if the endomorphism $\varphi_0=(f_0,g_0,z)$ of $K\langle x,y,z\rangle$
is a $z$-wild automorphism,
where $f_0,g_0$ are the linear in $x,y$ components of $f,g$, respectively.
This implies that the Anick automorphism is wild.
Recently Drensky and Yu \cite{DY4, DY5}
established the wildness of a big class of automorphisms and coordinates of
$K\langle x,y,z\rangle$. Many of them cannot be handled with direct application
of the methods of \cite{DY3} and \cite{U2}.
These results motivate the needs of systematic study of $z$-automorphisms
of $K\langle x,y,z\rangle$. As in the case of $z$-automorphisms of $K[x,y,z]$, they are
simpler than the arbitrary automorphisms of $K\langle x,y,z\rangle$ and
provide important examples and conjectures for
$\text{Aut }K\langle x,y,z\rangle$.

In the present paper we describe the structure of the group of
$z$-tame automorphisms of $K\langle x,y,z\rangle$ as the free
product of the groups of $z$-affine automorphisms and
$z$-triangular automorphisms amalgamating the intersection. We
also give algorithms which recognize $z$-tame automorphisms and
coordinates of $K\langle x,y,z\rangle$.
As an application, we show that all the $z$-automorphisms of the form
$\sigma_h=(x+zh(xz-zy,z),y+h(xz-zy,z)z)$ are $z$-wild when the polynomials
$h(xz-zy,z)$ are of positive degree
in $x$. This kind of automorphisms appear in \cite{DY4, DY5}
but the considerations there do not cover the case when $h(xz-zy,z)$
belongs to the square of the commutator ideal of $K\langle x,y,z\rangle$.
Besides, the polynomial $x+zh(xz-zy,z)$ is a $z$-wild coordinate.
Finally, we show that the $z$-endomorphisms of the form $\varphi=(x+u(x,y,z),y+v(x,y,z))$,
where $(u,v)\not=(0,0)$ and all monomials of $u$ and $v$ depend on both $x$ and $y$, are not automorphisms.
A partial case of this result was an essential
step in the proof of the theorem of Czerniakiewicz and Makar-Limanov
for the tameness of $\text{Aut }K\langle x,y\rangle$.
The paper may be considered as a continuation of our paper \cite{DY3}.

\section{The group of $z$-tame automorphisms}

We fix the field $K$ and consider the free associative algebra
$K\langle x,y,z\rangle$ in three variables.
We call the automorphism $\varphi$ of $K\langle x,y,z\rangle$ a $z$-automorphism if
$\varphi(z)=z$, and denote the automorphism group of the
$z$-automorphisms by $\text{Aut}_z\langle x,y,z\rangle$. Since we want to emphasize that we
work with $z$-automorphisms, we shall write $\varphi=(f,g)$,
omitting the third coordinate $z$.
The multiplication will be from right to left. If
$\varphi,\psi\in\text{Aut}_zK\langle x,y,z\rangle$, then in $\varphi\psi$ we
first apply $\psi$ and then $\varphi$. Hence, if $\varphi=(f,g)$
and $\psi=(u,v)$, then
\[
\varphi\psi=(u(f,g,z),v(f,g,z)).
\]
The $z$-affine and $z$-triangular
automorphisms of $K\langle x,y,z\rangle$ are, respectively, of the form
\[
\psi=(\alpha_{11}x+\alpha_{21}y+\alpha_{31}z+\beta_1,
\alpha_{12}x+\alpha_{22}y+\alpha_{32}z+\beta_2),
\]
$\alpha_{ij},\beta_j\in K$, the $2\times 2$ matrix $(\alpha_{ij})_{i,j=1,2}$
being invertible,
\[
\rho=(\alpha_1x+p_1(y,z),\alpha_2y+p_2(z)),
\]
$\alpha_j\in K^{\ast}$, $p_1\in K\langle y,z\rangle$, $p_2\in K[z]$.
The affine and the triangular $z$-automorphisms generate,
respectively, the subgroups $A_z$ and $B_z$ of $\text{Aut}_zK\langle x,y,z\rangle$.
We denote by $\text{TAut}_zK\langle x,y,z\rangle$ the group of $z$-tame
automorphisms which is generated by the $z$-affine and $z$-triangular
automorphisms. Of course, we may define the $z$-affine
automorphisms as the $z$-automorphisms of the form $\psi=(f,g)$, where
the polynomials $f,g\in K\langle x,y,z\rangle$ are linear in $x$ and $y$.
But, as we commented in \cite{DY3}, this definition is
not convenient. For example, the Anick
automorphism is affine in this sense but is wild.

In the commutative case, the $z$-automorphisms of $K[x,y,z]$ are simply the
automorphisms of the $K[z]$-algebra $K[z][x,y]$.
A result of Wright \cite{Wr} states that
over any field $K$ the group $\text{\rm TAut}_zK[x,y,z]$
has the amalgamated free product structure
\[
\text{\rm TAut}_zK[x,y,z]=A_z\ast_{C_z}B_z,
\]
where $A_z$ and $B_z$ are defined as in the case of $K\langle x,y,z\rangle$
and $C_z=A_z\cap B_z$.
(The original statement in \cite{Wr} holds in a more general situation.
In the case of $K[x,y,z]$ it involves affine and linear automorphisms
with coefficients from $K[z]$ but this is not
essential because every invertible matrix with entries in $K[z]$ is a product
of elementary and diagonal matrices.)

Every $z$-tame automorphism $\varphi$ of $K\langle x,y,z\rangle$
can be presented as a product in the form
(\ref{canonical form of automorphisms})
where $\psi_i\in A_z$, $\rho_i\in B_z$ ($\varepsilon_1$ and $\varepsilon_m$ are
equal to 0 or 1), and, if $\varphi$ does not belong to the union of $A_z$ and $B_z$,
we may assume that $\psi_i\in A_z\backslash B_z$, $\rho_i\in B_z\backslash A_z$.
Fixing the linear nontriangular $z$-automorphism $\tau=(y,x)$,
we can present $\varphi$ in the canonical form
\begin{equation}\label{simplified canonical form of automorphisms}
\varphi=\rho_n\tau\cdots\tau\rho_1\tau\rho_0,
\end{equation}
where $\rho_0,\rho_1,\ldots,\rho_n\in B_z$ and only $\rho_0$ and $\rho_n$
are allowed to belong to $A_z$, see for example p. 350 in \cite{C2}.
Let
\[
\rho_i=(\alpha_ix+p_i(y,z),\beta_iy+r_i(z)),\quad
\alpha_i,\beta_i\in K^{\ast},p_i\in K\langle y,z\rangle,
r_i\in K[z].
\]
Using the equalities for compositions of automorphisms
\[
(\alpha x+p(y,z),\beta y+r(z))=
(x+\alpha^{-1}(p(y,z)-p(0,z)),y)(\alpha x+p(0,z),\beta y+r(z)),
\]
\[
(\alpha x+p(z),\beta y+r(z))\tau=(\beta y+r(z),\alpha x+p(z))
=\tau(\beta x+r(z),\alpha y+p(z)),
\]
$p(z),r(z)\in K[z]$,
we can do further simplifications in (\ref{simplified canonical form of automorphisms}),
assuming that $\rho_1,\ldots,\rho_{n-1}$ are not affine and, together with $\rho_n$,
are of the form $\rho_i=(x+p_i(y,z),y)$ with $p_i(0,z)=0$ for all $i=1,\ldots,n$.
We also assume that $\rho_0=(\alpha_0x+p_0(y,z),\beta_0y+r_0(z))$.
The condition that $\rho_1,\ldots,\rho_{n-1}$ are not affine means that
$\text{deg}_yp_i(y,z)\geq 1$ and if $\text{deg}_yp_i(y,z)=1$, then
$\text{deg}_zp_i(y,z)\geq 1$, $i=1,\ldots,n-1$.

The following result shows that the structure of the group
of $z$-tame automorphisms of $K\langle x,y,z\rangle$ is similar to the structure of
the group of $z$-tame automorphsims fo $K[x,y,z]$.

\begin{theorem}\label{structure of z-TAut}
Over an arbitrary field $K$,
the group $\text{\rm TAut}_zK\langle x,y,z\rangle$ of $z$-tame
automorphisms of $K\langle x,y,z\rangle$
is isomorphic to the free product $A_z\ast_{C_z}B_z$ of the group $A_z$
of the $z$-affine automorphisms and the group $B_z$ of $z$-triangular automorphisms
amalgamating their intersection $C_z=A_z\cap B_z$.
\end{theorem}

\begin{proof}
We define a bidegree of $K\langle x,y,z\rangle$ assuming that the monomial $w$ is
of bidegree $\text{bideg }w=(d,e)$ if
$\text{deg}_xw+\text{deg}_yw=d$ and $\text{deg}_zw=e$.
We order the bidegrees $(d,e)$ lexicographically, i.e.,
$(d_1,e_1)>(d_2,e_2)$ means that either $d_1>d_2$ or
$d_1=d_2$ and $e_1>e_2$. We denote by $\overline{p}$ the leading
bihomogeneous component  of the nonzero polynomial $p(x,y,z)$.
Let $\varphi=(f,g)$ be in the form
(\ref{simplified canonical form of automorphisms}), with all the
restrictions fixed above, and let $q_i(y,z)$ be the leading
component of $p_i(y,z)$. Direct computations give that, if
$\rho_n$ is not linear and $p_0(y,z)\not=\gamma_0y+p'_0(z)$ in
$\rho_0=(\alpha_0x+p_0(y,z),\beta_0y+r_0(z))$, then
\begin{equation}\label{the leading terms of f and g}
\begin{array}{c}
\overline{f}=q_0(q_1(\ldots q_{n-1}(q_n(y,z),z)\ldots,z),z),\\
\\
\overline{g}=\beta_0q_1(\ldots q_{n-1}(q_n(y,z),z)\ldots,z),
\end{array}
\end{equation}
and $\text{bideg }\overline{f}>(1,0)$. Hence $\varphi$ is not the
identity automorphism. Similar considerations work when at least one of
the automorphisms $\rho_0$ and $\rho_n$ is affine. For example, if
$\rho_0=(\alpha_0+\gamma_0y+p'_0(z),\beta_0y+r_0(z))$,
$\gamma_0\in K^{\ast}$, and $\text{bideg }p_n(y,z)>(1,0)$, then
\[
\overline{f}=\gamma_0q_1(\ldots q_{n-1}(q_n(y,z),z)\ldots,z),
\]
\[
\overline{g}=\beta_0q_1(\ldots q_{n-1}(q_n(y,z),z)\ldots,z).
\]
If $\text{bideg }p_0(y,z)>(1,0)$ and $\rho_n=(x+\gamma_ny,y)$,
$\gamma_n\in K^{\ast}$, then
\[
\overline{f}=q_0(q_1(\ldots q_{n-1}(q_n(x+\gamma y,z),z)\ldots,z),z),
\]
\[
\overline{g}=\beta_0q_1(\ldots q_{n-1}(q_n(x+\gamma y,z),z)\ldots,z).
\]
In all the cases, $\varphi$ is not the identity automorphism.
Hence, if $\varphi$ has a nontrivial presentation in the form
(\ref{simplified canonical form of automorphisms}), then it is
different from the identity automorphism, and we conclude that
$\text{TAut}_zK\langle x,y,z\rangle$ is a free product with amalgamation
of the groups $A_z$ and $B_z$.
\end{proof}

Following our paper \cite{DY3} we identify the group
of $z$-automorphisms which are linear in $x$ and $y$
with the group $GL_2(K[z_1,z_2])$.
Let $f\in K\langle x,y,z\rangle$ be linear in $x,y$.
Then $f$ has the form
\[
f=\sum\alpha_{ij}z^ixz^j+\sum\beta_{ij}z^iyz^j,\quad \alpha_{ij},\beta_{ij}\in K.
\]
The $z$-derivatives
$f_x$ and $f_y$ are defined by
\[
f_x=\sum\alpha_{ij}z_1^iz_2^j,\quad
f_y=\sum\beta_{ij}z_1^iz_2^j.
\]
Here $f_x$ and $f_y$ are in $K[z_1,z_2]$ and are polynomials in two commuting variables.
The $z$-Jacobian matrix of the linear $z$-endomorphism $\varphi=(f,g)$
of $K\langle x,y,z\rangle$ is defined as
\[
J_z(\varphi)=\left(
\begin{matrix}
f_x&g_x\\
f_y&g_y\\
\end{matrix}
\right).
\]
By \cite{DY3} the mapping $\varphi\to J_z(\varphi)$
is an isomorphism of the group of the $z$-automorphisms which
are linear in $x,y$ and $GL_2(K[z_1,z_2])$. Also, such an
automorphism is $z$-tame if and only if
its $z$-Jacobian matrix belongs to $GE_2(K[z_1,z_2])$.
(By the further development of this result by Umirbaev \cite{U2},
the $z$-wild automorphisms of the considered type are wild
also as automorphisms of $K\langle x,y,z\rangle$.)

\begin{corollary}
The group $\text{\rm TAut}_zK\langle x,y,z\rangle$
is isomorphic to the free product with amalgamation
$GE_2(K[z_1,z_2])\ast_{C_1}B_z$, where
$GE_2(K[z_1,z_2])$ is identified as above with the group of $z$-tame automorphisms
which are linear in $x$ and $y$, and
$C_1=GE_2(K[z_1,z_2])\cap B_z$.
\end{corollary}

\begin{proof}
Everything follows from the observations that: (i) in
the form (\ref{simplified canonical form of automorphisms}),
$\rho_j\tau\cdots\tau\rho_i\in GE_2(K[z_1,z_2])$ if and only if
all $\rho_j,\ldots,\rho_i$ belong to $GE_2(K[z_1,z_2])$;
(ii) $\rho_j\tau\cdots\tau\rho_i\in C_1$ if and only if $i=j$ and
$\rho_i\in GE_2[z_1,z_2]$; (iii) $\tau\in GE_2(K[z_1,z_2])$.
\end{proof}

\section{Recognizing $z$-tame automorphisms and coordinates}

Now we use Theorem \ref{structure of z-TAut} to present algorithms
which recognize $z$-tame automorphisms and
coordinates of $K\langle x,y,z\rangle$. Of course, in all
algorithms we assume that the field $K$ is constructive.
We start with an algorithm
which determines whether a $z$-endomorphism of $K\langle x,y,z\rangle$
is a $z$-tame automorphism. The main idea is similar to that of
the well known algorithm which decides whether an endomorphism of
$K[x,y]$ is an automorphism, see Theorem 6.8.5 in \cite{C2}, but the realization
is more sophisticated. In order to simplify the considerations, we
shall use the trick introduced by Formanek \cite{Fo} in his
construction of central polynomials of matrices.

Let $H_n$ be the subspace of $K\langle x,y,z\rangle$ consisting of all
polynomials which are homogeneous of degree $n$ with respect to $x$ and $y$.
We define an action of $K[t_0,t_1,\ldots,t_n]$ on $H_n$ in the following way. If
\[
w=z^{a_0}u_1z^{a_1}u_2\cdots z^{a_{n-1}}u_nz^{a_n},
\]
where $u_i=x$ or $u_i=y$, $i=1,\ldots,n$, then
\[
t_0^{b_0}t_1^{b_1}\cdots t_n^{b_n}\ast w=
z^{a_0+b_0}u_1z^{a_1+b_1}u_2\cdots z^{a_{n-1}+b_{n-1}}u_nz^{a_n+b_n},
\]
and then extend this action by linearity. Clearly, $H_n$ is a free
$K[t_0,t_1,\ldots,t_n]$-module with basis consisting of the $2^n$
monomials $u_1\cdots u_n$, where $u_i=x$ or $u_i=y$.
The proof of the following lemma is obtained by easy direct computation.

\begin{lemma}\label{trick of Formanek}
Let $\beta\in K^{\ast}$,
\begin{equation}\label{this is v}
v(x,y,z)=\sum\theta_i(t_0,t_1,\ldots,t_k)\ast u_{i_1}\cdots u_{i_k}\in H_k,
\end{equation}
\begin{equation}\label{this is q}
q(y,z)=\omega(t_0,t_1,\ldots,t_d)\ast y^d\in H_d,
\end{equation}
where $\theta_i\in K[t_0,t_1,\ldots,t_k]$, $\omega\in K[t_0,t_1,\ldots,t_d]$,
$u_{i_j}=x$ or $u_{i_j}=y$. Then
\begin{equation}\label{this is u}
\begin{array}{c}
u(x,y,z)=q(v(x,y,z)/\beta,z)
=\omega(t_0,t_d,t_{2d},\ldots,t_{kd})/\beta^d\\
\\
\left(\sum\theta_i(t_0,t_1,\ldots,t_k)\ast u_{i_1}\cdots u_{i_k}\right)\\
\\
\left(\sum\theta_i(t_k,t_{k+1},\ldots,t_{2k})\ast u_{i_1}\cdots u_{i_k}\right)\cdots\\
\\
\left(\sum\theta_i(t_{k(d-1)},t_{k(d-1)+1},\ldots,t_{kd})\ast u_{i_1}\cdots u_{i_k}\right).\\
\end{array}
\end{equation}
\end{lemma}

\begin{algorithm}\label{algorithm for z-tame automorphisms}
Let $\varphi=(f,g)$ be a $z$-endomorphism of $K\langle x,y,z\rangle$.
We make use of the bidegree defined in the proof of Theorem
\ref{structure of z-TAut}.

{\it Step 0}. If some of the polynomials $f,g$ depends on $z$ only,
then $\varphi$ is not an automorphism.

{\it Step 1}. Let $u,v$ be the homogeneous components of highest bidegree of $f,g$,
respectively. If both $u,v$ are of bidegree $(1,0)$, i.e., linear, then we check
whether they are linearly
independent. If yes, then $\varphi$ is a product of a linear automorphism
(from $GL_2(K)$) and a translation $(x+p(z),y+r(z))$. If
$u,v$ are linearly dependent, then $\varphi$ is not an automorphism.

{\it Step 2}. Let $\text{bideg }u>(1,0)$ and
$\text{bideg }u\geq \text{bideg }v$. Hence $u\in H_l$, $v\in H_k$ for some $k$ and $l$.
Taking into account (\ref{the leading terms of f and g}),
we have to check whether $l=kd$ for a positive integer $d$
and to decide whether $u=q(v/\beta,z)$ for some
$\beta\in K^{\ast}$ and some $q(y,z)\in H_d$. In the notation of
Lemma \ref{trick of Formanek}, we know
$u$ in (\ref{this is u}) and $v$ in (\ref{this is v}) up to the multiplicative constant $\beta$.
Hence, up to $\beta$, we know the polynomials
$\theta_i(t_0,t_1,\ldots,t_n)$ in the presentation of $v$.
We compare some of the nonzero polynomial coefficients of
$u=\sum\lambda_j(t_0,\ldots,t_{kd})u_{j_1}\cdots u_{i_{kd}}$
with the corresponding coefficient of $q(v/\beta,z)$.
Lemma \ref{trick of Formanek} allows to find explicitly, up to
the value of $\beta^d$, the polynomial
$\omega(t_0,t_1,\ldots,\omega_d)$ in (\ref{this is q}) using the usual division of polynomials.
If $l=kd$ and $u=q(v/\beta,z)$, then we replace
$\varphi=(f,g)$ with $\varphi_1=(f-q(g/\beta,z),g)$.
Then we apply Step 0 to $\varphi_1$.
If $u$ cannot be presented in the desired form,
then $\varphi$ is not an automorphism.

{\it Step 3}. If $\text{bideg }v>(1,0)$ and
$\text{bideg }u<\text{bideg }v$, we have similar considerations, as
in Step 2, replacing
$\varphi=(f,g)$ with $\varphi_1=(f,g-q(f/\alpha,z))$ for suitable
$q(y,z)$.
Then we apply Step 0 to $\varphi_1$.
If $v$ cannot be presented in this form,
then $\varphi$ is not an automorphism.
\end{algorithm}

\begin{corollary}\label{nonmetabelian generalization of Anick}
Let $h(t,z)\in K\langle t,z\rangle$ and let $\text{\rm deg}_uh(u,z)>0$.
Then
\[
\sigma_h=(x+zh(xz-zy,z),y+h(xz-zy,z)z,z)
\]
is a $z$-wild automorphism of $K\langle x,y,z\rangle$.
\end{corollary}

\begin{proof}
It is easy to see that $\sigma_h$ is a $z$-automorphism
of $K\langle x,y,z\rangle$ with inverse $\sigma_{-h}$.
We apply Algorithm \ref{algorithm for z-tame automorphisms}. Let $w$ be the
homogeneous component of highest bidegree of $h(xz-zy,z)$. Clearly,
$w$ has the form $w=\overline{h}(xz-zy,z)=q(xz-zy,z)$ for some bihomogeneous
polynomial $q(t,z)\in K\langle t,z\rangle$. The leading components
of the coordinates of $\sigma_h$ are $zq(xz-zy,z)$ and $q(xz-zy,z)z$, and are of
the same bidegree. If $\sigma_h$ is a $z$-tame automorphism, then we can
reduce the bidegree using a linear transformation, which is impossible because
$zq(xz-zy,z)$ and $q(xz-zy,z)z$ are linearly independent.
\end{proof}

The algorithm in Theorem 6.8.5 in \cite{C2} which recognizes the automorphisms of
$K[x,y]$ can be easily modified to recognize the coordinates of $K[x,y]$.
Such an algorithm is explicitly stated
in \cite{SY3}, where Shpilrain and Yu established an algorithm
which gives a canonical form, up to automorphic equivalence,
of a class of polynomials in $K[x,y]$. (The automorphic equivalence problem
for $K[x,y]$ asks how to decide whether, for two given polynomials $p,q\in K[x,y]$,
there exists an automorphism $\varphi$ such that $q=\varphi(p)$. It was solved
over $\mathbb C$ by Wightwick \cite{Wi} and, over an arbitrary algebraically closed
constructive field $K$, by Makar-Limanov, Shpilrain, and Yu \cite{MLSY}.)
When $\text{char }K=0$,
Shpilrain and Yu \cite{SY1} gave a very simple algorithm
which decides whether a polynomial
$f(x,y)\in K[x,y]$ is a coordinate. Their approach is based on an idea of
Wright \cite{Wr} and the Euclidean division algorithm applied for
the partial derivatives of a polynomial in $K[x,y]$.
Using the isomorphism of $\text{Aut }K[x,y]$ and $\text{Aut }K\langle x,y\rangle$
and reducing the considerations to the case of $K[x,y]$,
Shpilrain and Yu \cite{SY2}
found the first algorithm which recognizes the coordinates
of $K\langle x,y\rangle$. Now we want to modify
Algorithm \ref{algorithm for z-tame automorphisms}
to decide whether a polynomial $f(x,y,z)$ is a $z$-tame coordinate
of $K\langle x,y,z\rangle$.

Note, that if $\varphi=(f,g)$ and $\varphi'=(f,g')$ are two $z$-automorphisms
of $K\langle x,y,z\rangle$ with the
same first coordinate, then $\varphi^{-1}\varphi'$ fixes $x$. Hence
$\varphi^{-1}\varphi'=(x,g'')$ and, obligatorily, $g''=\beta y+r(x,z)$.
In this way, if we know one $z$-coordinate mate $g$ of $f$, then we are able to find
all other $z$-coordinate mates.
These arguments and Corollary \ref{nonmetabelian generalization of Anick} give immediately:

\begin{corollary}\label{coordinates of nonmetabelian generalization of Anick}
Let $h(t,z)\in K\langle t,z\rangle$ and let $\text{\rm deg}_uh(u,z)>0$.
Then $f(x,y,z)=x+zh(xz-zy,z)$
is a $z$-wild coordinate of $K\langle x,y,z\rangle$.
\end{corollary}

\begin{theorem}\label{algorithm for z-tame coordinates}
There is an algorithm which decides whether a polynomial
$f(x,y,z)\in K\langle x,y,z\rangle$ is a $z$-tame coordinate.
\end{theorem}

\begin{proof}
We start with the analysis of the behavior of the first coordinate $f$
of $\varphi$ in (\ref{simplified canonical form of automorphisms}).
Let $h$ be the first coordinate of
$\psi=\rho_{n-1}\tau\cdots\tau\rho_1\tau\rho_0$ and let, as in
(\ref{simplified canonical form of automorphisms}), $\rho_n=(x+p_n(y,z),y)$
and $p_n(0,z)=0$. Then
\begin{equation}\label{relation between f and h}
f(x,y,z)=\rho_n\tau(h(x,y,z))=h(y,x+p_n(y,z),z).
\end{equation}
In order to make the inductive step, we have to recover the polynomials $h(x,y,z)$ and $p_n(y,z)$
or, at least their leading components with respect to a suitable grading.

For a pair of positive integers
$(a,b)$, we define the $(a,b)$-bidegree of a monomial
$w\in K\langle x,y,z\rangle$ by
\[
\text{bideg}_{(a,b)}w=(a\text{deg}_xw+b\text{deg}_yw,\text{deg}_zw)
\]
and order the bidegrees in the lexicographic order, as in
Algorithm \ref{algorithm for z-tame automorphisms}. For a nonzero polynomial
$f\in K\langle x,y,z\rangle$ we denote by
$\vert f\vert_{(a,b)}$ the homogeneous component of maximal $(a,b)$-bidegree.
We write $\varphi=(f,g)\in \text{TAut}_zK\langle x,y,z\rangle$ in the
form (\ref{simplified canonical form of automorphisms}).
Let us assume again that $\text{bideg }p_i(y)>(1,0)$ for all $i=0,1,\ldots,n$, and
let $h$ be the first coordinate of
$\psi=\rho_{n-1}\tau\cdots\tau\rho_1\tau\rho_0$. Then the highest
bihomogeneous component of $h$ is
\[
\overline{h}(y,z)=q_0(q_1(\ldots (q_{n-1}(y,z),z)\ldots),z).
\]
The homogeneous component of maximal $(d_n,1)$-bidegree of $x+p_n(y,z)$
is $\vert x+q_n(y,z)\vert_{(d_n,1)}=x+\xi_ny^{d_n}$ if $\text{deg}_zq_n(y,z)=0$ and
$\vert x+q_n(y,z)\vert_{(d_n,1)}=q_n(y,z)$ if $\text{deg}_zq_n(y,z)>0$.
Direct calculations give
\[
\vert f\vert_{(d_n,1)}=\vert \rho_n\tau(\overline{h})\vert_{(d_n,1)}
=\vert\overline{h}(x+q_n(y,z)\vert_{(d_n,1)}.
\]
If $f'(x,z)$ and $f''(y,z)$ are the components of $f(x,y,z)$ which do not depend on $y$ and $x$,
respectively, we can recover the degree $d_n$ of $p_n(y,z)$ as the quotient
$d_n=\text{deg}_xf'/\text{deg}_yf''$. Now the problem is to recover $q_n(y,z)$ and $\overline{h}(y,z)$.
Since $\overline{h}(y,z)$ does not depend on $x$, we have that
\[
\overline{h}(y,z)=\overline{h(x,y,z)}=\overline{h(0,y,z)}.
\]
From the equality (\ref{relation between f and h}) and the condition $p_n(0,z)=0$
we obtain that
\[
f(x,0,z)=h(0,x+p_n(0,z),z)=h(0,x,z).
\]
Hence $h(0,y,z)=f(y,0,z)$ and we are able to find $\overline{h}(y,z)$.
We write $\overline{h}$ and $\overline{q_n}$ in the form
\[
\overline{h}(y,z)=\theta(t_0,t_1,\ldots,t_k)\ast y^k,\quad
q_n(y,z)=\omega(t_0,t_1,\ldots,t_d)\ast y^d,
\]
where $\theta(t_0,t_1,\ldots,t_k)\in K[t_0,t_1,\ldots,t_k]$ is known explicitly and
$\omega(t_0,t_1,\ldots,t_d)\in K[t_0,t_1,\ldots,t_d]$.
Similarly, the part of the component of maximal bidegree of $f(x,y,z)$ which does not depend on $x$
has the form
\[
\overline{f''}(y,z)=\zeta(t_0,t_1,\ldots,t_{kd})\ast y^{kd},\quad
\zeta(t_0,t_1,\ldots,t_{kd})\in K[t_0,t_1,\ldots,t_{kd}].
\]
Since $\overline{h}(q_n(y,z),z)=\overline{f''}(y,z)$, by Lemma \ref{trick of Formanek} we obtain
\[
\zeta(t_0,t_1,\ldots,t_{kd})=\theta(t_0,t_d,t_{2d},\ldots,t_{kd})
\omega(t_0,t_1,\ldots,t_d)
\]
\[
\omega(t_d,t_{d+1},\ldots,t_{2d})\cdots
\omega(t_{(k-1)d},t_{(k-1)d+1},\ldots,t_{kd}).
\]
Here we know $\zeta$ and $\theta$ and want to determine $\omega$.
Let
\[
\zeta'(t_0,t_1,\ldots,t_{kd})=\zeta(t_0,t_1,\ldots,t_{kd})/\theta(t_0,t_d,t_{2d},\ldots,t_{kd})
\]
\[
=\omega(t_0,t_1,\ldots,t_d)
\omega(t_d,t_{d+1},\ldots,t_{2d})\cdots
\omega(t_{(k-1)d},t_{(k-1)d+1},\ldots,t_{kd}).
\]
The greatest common divisor of the polynomials
$\zeta'(t_0,t_1,\ldots,t_{kd})$ and $\zeta'(t_{(k-1)d},t_{(k-1)d+1},\ldots,t_{(2k-1)d})$
in $K[t_0,t_1,\ldots,t_{(2k-1)d}]$ is equal, up to a multiplicative constant $\beta$,
to $\omega(t_{(k-1)d},t_{(k-1)d+1},\ldots,t_{kd})$.
Hence the knowledge of $\zeta'$ allows to determine $\beta\omega(t_0,t_1,\ldots,t_k)$
as well as the value of $\beta^d$. This means that we know also all the possible values of $\beta$
and the polynomial $q_n(y,z)$.
Now we apply on $f(x,y,z)$ the $z$-automorphism $\sigma=(x-q_n(y,z),y)$. Since
$f(x,y,z)-\overline{h}(x+q_n(y,z),z)$ is lower in the $(d_n,1)$-biordering
than $f(x,y,z)$ itself, we may replace $f$ with
$\sigma(f)$ and to make the next step.
The considerations are almost the same when
some of the automorphisms $\rho_0$ and $\rho_n$ is affine.
For example, if $f=\varphi(x)$ and
$\rho_n=(x+\gamma y,y)$, $\gamma\in K$, in (\ref{simplified canonical form of automorphisms}),
then the leading bihomogeneous component of $h=\tau\rho_n^{-1}(f)$ does not depend on $y$,
and we can do the next step. If $f$ is a $z$-tame coordinate, then
the above process will stop when we reduce $f$ to a polynomial in the form
$\alpha x+p(y,z)$. If $f$ is not a $z$-tame coordinate, then the process will also stop by different reason.
In some step we shall reduce $f(x,y,z)$ to a polynomial $f_1(x,y,z)$.
It may turn out that the degree $d=\text{deg}_xf_1(x,0,z)/\text{deg}_yf_1(0,y,z)$
is not integer. Or, the commutative polynomials $\theta$ and $\omega$ corresponding to $f_1$ do not exist.
\end{proof}

The following corollary is stronger than Corollary \ref{coordinates of nonmetabelian generalization of Anick}.

\begin{corollary}\label{stronger nonmetabelian generalization of Anick}
Let $h(t,z)\in K\langle t,z\rangle$ and let $\text{\rm deg}_uh(u,z)>0$.
Then $f(x,y,z)=x+h(xz-zy,z)$
is not a $z$-tame coordinate of $K\langle x,y,z\rangle$.
\end{corollary}

\begin{proof}
We apply the algorithm in the proof of Theorem \ref{algorithm for z-tame coordinates}.
Let $f(x,y,z)$ be a $z$-tame coordinate and let $h'(x,z)=h(xz,z)$ and $h''(y,z)=h(-zy,z)$
be the polynomials obtained from $h(xz-zy,z)$
replacing, respectively, $y$ and $x$ by 0. Clearly, $\text{bideg}_xh'=\text{bideg}_yh''$. Hence,
as in the proof of Theorem \ref{algorithm for z-tame coordinates} we can replace
$f(x,y,z)$ with $\sigma(f)$, where $\sigma=(x-\alpha y,y)$,
for a suitable $\alpha\in K^{\ast}$, and the leading bihomogeneous component of
$\sigma(f)$ in the $(1,1)$-ordering does not depend on $y$.
But this brings to a contradiction. If $h_1(t,z)\in K\langle t,z\rangle$
is homogeneous with respect to $t$, and
\[
h_1((x-\gamma y)z-zy,z)=h_2(x,z)
\]
for some $h_2(x,z)$, then, replacing $x$ with $0$, we obtain
$h_1(-(\gamma yz+zy),z)=0$, which is impossible.
\end{proof}

\begin{remark}
In Corollary \ref{stronger nonmetabelian generalization of Anick}, we cannot guarantee that
the polynomial $f(x,y,z)=x+h(xz-zy,z)$ is a $z$-coordinate at all.
For example, let $f(x,y,z)=x+(xz-zy)$ be a $z$-coordinate
with a coordinate mate $g(x,y,z)$. If $g_1(x,y,z)$ is the linear in $x,y$ component of $g$, then
$\varphi_1=(f,g_1)$ is also a $z$-automorphism. Then, for suitable polynomials
$c, d\in K[z_1,z_2]$,
the matrix
\[
J_z(\varphi_1)=\left(
\begin{matrix}
1+z_2&c(z_1,z_2)\\
-z_1&d(z_1,z_2)\\
\end{matrix}
\right)
\]
is invertible. If we replace $z_1$ with $0$
in its determinant $\text{det}(J_z)=(1+z_2)d(z_1,z_2)-z_1c(z_1,z_2)$ we obtain
that $(1+z_2)d_2(0,z_2)\in K^{\ast}$ which is impossible.
\end{remark}

\section{Endomorphisms which are not automorphisms}

In this section we shall establish a $z$-analogue of the following proposition
which is the main step of the proof of the
theorem of Czerniakiewicz \cite{Cz} and Makar-Limanov \cite{ML1, ML2}
for the tameness of the automorphisms of $K\langle x,y\rangle$.

\begin{proposition}\label{main step of the proof of tameness}
Let $\varphi=(x+u,y+v)$ be an endomorphism of $K\langle x,y\rangle$,
where $u,v$ are in the commutator ideal of $K\langle x,y\rangle$
and at least one of them is different from 0.
Then $\varphi$ is not an automorphism of $K\langle x,y\rangle$.
\end{proposition}

An essential moment in its proof,
see the book by Cohn \cite{C2}, is the following lemma.

\begin{lemma}\label{lemma for bihomogenoeus polynomials}
If $f,g\in K\langle x,y\rangle$ are two bihomogeneous polynomials, then they either generate
a free subalgebra of $K\langle x,y\rangle$ or, up to multiplicative constants, both are
powers of the same bihomogeneous element of $K\langle x,y\rangle$.
\end{lemma}

We shall prove a weaker version of the lemma for
$K\langle x,y,z\rangle$ which will be sufficient
for our purposes.

\begin{lemma}\label{analogue of the lemma from the book by Cohn}
Let $(0,0)\not=(a,b)\in {\mathbb Z}^2$ and let $f_1,f_2\in K\langle x,y,z\rangle$ be bihomogeneous
with respect to the $(a,b)$-degree of $K\langle x,y,z\rangle$, i.e.,
$a\text{\rm deg}_xw+b\text{\rm deg}_yw$ is the same for all monomials of $f_1$, and similarly
for $f_2$. If $f_1$ and $f_2$ are algebraically dependent, then both
$\text{\rm deg}_{(a,b)}f_1$ and $\text{\rm deg}_{(a,b)}f_2$ are either nonnegative or nonpositive.
\end{lemma}

\begin{proof}
Let $v(f_1,f_2,z)=0$ for some nonzero polynomial $v(u_1,u_2,z)\in K\langle u_1,u_2,z\rangle$.
We may assume that both $f_1,f_2$ depend not on $z$ only.
We fix a term-ordering on $K\langle x,y,z\rangle$. Let $\tilde f_1$ and $\tilde f_2$ be the leading
monomials of $f_1$ and $f_2$, respectively. For each monomial
$z^{k_0}u_{i_1}z^{k_1}\cdots z^{k_{s-1}}u_{i_s}z^{k_s}\in K\langle u_1,u_2,z\rangle$
the leading monomial of $z^{k_0}f_{i_1}z^{k_1}\cdots z^{k_{s-1}}f_{i_s}z^{k_s}\in K\langle x,y,z\rangle$
is $z^{k_0}\tilde f_{i_1}z^{k_1}\cdots z^{k_{s-1}}\tilde f_{i_s}z^{k_s}$. Hence, the algebraic dependence
of $f_1$ and $f_2$ implies that two different monomials
$z^{k_0}\tilde f_{i_1}z^{k_1}\cdots z^{k_{s-1}}\tilde f_{i_s}z^{k_s}$
and $z^{l_0}\tilde f_{j_1}z^{l_1}\cdots z^{l_{t-1}}\tilde f_{j_t}z^{j_t}$ are equal.
We write $\tilde f_1=z^{p_1}g_1z^{q_1}$ and $\tilde f_2=z^{p_2}g_2z^{q_2}$, where $g_1, g_2$ do not start
and do not end with $z$. After some cancelation in the equation
\[
z^{k_0}\tilde f_{i_1}z^{k_1}\cdots z^{k_{s-1}}\tilde f_{i_s}z^{k_s}=
z^{l_0}\tilde f_{j_1}z^{l_1}\cdots z^{l_{t-1}}\tilde f_{j_t}z^{j_t}
\]
we obtain a relation of the form
\begin{equation}\label{two monomials are equal}
g_{a_1}z^{m_1}\cdots z^{m_{k-1}}g_{a_k}z^{m_k}=
g_{b_1}z^{n_1}\cdots z^{n_{l-1}}g_{b_l}z^{n_l},
\end{equation}
with different $g_{a_1}$ and $g_{b_1}$. Hence, if $\text{deg }g_1\geq \text{deg }g_2$, then
$g_1=g_2g_3$ for some monomial $g_3$ (and $g_2=g_1g_3$ if $\text{deg }g_1< \text{deg }g_2$).
Again, $g_2$ and $g_3$ satisfy a relation of the form (\ref{two monomials are equal}).
Since $\text{deg }g_1\geq \text{deg }g_2>0$, we obtain
$\text{deg }g_1+\text{deg }g_2>\text{deg }g_1=\text{deg }g_2+\text{deg }g_3$. Applying inductive arguments,
we derive that both $\text{\rm deg}_{(a,b)}g_2$ and $\text{\rm deg}_{(a,b)}g_3$
are either nonnegative or nonpositive, and the same holds for $f_1$ and $f_2$ because
$g_1=g_2g_3$, $\text{\rm deg}_{(a,b)}g_1=\text{\rm deg}_{(a,b)}g_2+\text{\rm deg}_{(a,b)}g_3$, and
$\text{\rm deg}_{(a,b)}f_i=\text{\rm deg}_{(a,b)}g_i$, $i=1,2$.
\end{proof}

The condition that $u(x,y)$ and $v(x,y)$ belong to the commutator ideal of $K\langle x,y\rangle$, as
in Proposition \ref{main step of the proof of tameness}, immediately implies that
all monomials of $u$ and $v$ depend on both $x$ and $y$, as required in the
following theorem.

\begin{theorem}
The $z$-endomorphisms of the form
\[
\varphi=(x+u(x,y,z),y+v(x,y,z)),
\]
where $(u,v)\not=(0,0)$ and all monomials of $u$ and $v$ depend on both $x$ and $y$, are not automorphisms
of $K\langle x,y,z\rangle$.
\end{theorem}

\begin{proof}
The key moment in the proof of Proposition \ref{main step of the proof of tameness}
is the following. If $\varphi=(x+u,y+v)$ is an endomorphism of $K\langle x,y\rangle$,
where $u,v$ are in the commutator ideal of $K\langle x,y\rangle$
and at least one of them is different from 0, then there exist two integers $a$ and $b$ such that
$(a,b)\not=(0,0)$ and $a\leq 0\leq b$ with the property that $\text{deg}_{(a,b)}(x+u)=\text{deg}_{(a,b)}x=a$
and $\text{deg}_{(a,b)}(y+v)=\text{deg}_{(a,b)}y=b$. Ordering in a suitable way the $(a,b)$-bidegrees,
one concludes that the $(a,b)$-degrees of the leading bihomogeneous components
of $x+u$ and $y+v$ are with different signs. Then Lemma \ref{lemma for bihomogenoeus polynomials}
shows that these leading components are algebraically independent and bidegree arguments
as in the proof of Proposition \ref{main step of the proof of tameness} give that
$\varphi$ cannot be an automorphism. We repeat verbatim these arguments, working with the same
$(a,b)$-(bi)degree and bidegree ordering Proposition \ref{main step of the proof of tameness},
without counting the degree of $z$. In the final step, we use
Lemma \ref{analogue of the lemma from the book by Cohn}
instead of Lemma \ref{lemma for bihomogenoeus polynomials}.
\end{proof}

\end{document}